\documentclass{endm}
\usepackage{endmmacro}
\usepackage{graphicx}

\begin{document}

% DO NOT REMOVE: Creates space for Elsevier logo, ScienceDirect logo
% and ENDM logo
\begin{verbatim}\end{verbatim}\vspace{2.5cm}

\begin{frontmatter}

\title{FHCP Challenge Set: The First Set of Structurally Difficult Instances of the Hamiltonian Cycle Problem}

\author{Michael Haythorpe\thanksref{myemail}}
\address{College of Science and Engineering\\ Flinders University\\ Tonsley Park, Australia}

\thanks[myemail]{Email:
   \href{michael.haythorpe@flinders.edu.au} {\texttt{\normalshape
   michael.haythorpe@flinders.edu.au}}}

\begin{abstract}
The FHCP Challenge Set, comprising of 1001 instances of Hamiltonian cycle problem, is introduced. This set is the first to contain instances of Hamiltonian cycle problem for which the primary difficulty is the underlying graph structure, rather than simply size. A summary of the kinds of graphs contained in the FHCP Challenge Set is given. A discussion of the results of the FHCP Challenge, a year-long competition to solve all instances of the FHCP Challenge Set first announced at the 59th Annual Meeting for the Australian Mathematical Society, is also included.
\end{abstract}


\end{frontmatter}

\section{Introduction}

The {\em Hamiltonian cycle problem} (HCP) is a, now classical, problem in graph theory which can be stated simply as follows. Given a graph $G$ containing $n$ vertices, does $G$ possess any simple cycles of length $n$? Such a simple cycle traversing every vertex in the graph is called a {\em Hamiltonian cycle} (HC), and graphs which contain at least one HC are called {\em Hamiltonian}. Although it was first posed in the 1850s, there has been a modern resurgence of interest in HCP. This has been in part because HCP was among the first problems proved to be NP-complete \cite{karp}, and also because of its very close relationship with the, even more famous, traveling salesman problem (TSP). Indeed, the latter can be posed as the question of finding a Hamiltonian cycle of optimal length in a weighted graph. For a summary of results on TSP, the interested reader is referred to Lawler et al \cite{lawler} which also contains a chapter on HCP written by Chv\'{a}tal.

There are now several efficient heuristics which are capable of solving HCP, including specialised HCP algorithms \cite{chalaturnyk,SLH,eppstein,ali,gondzio,gazettepaper,pettersson,branchbound} as well as TSP algorithms which can be adapted to solve HCP \cite{concorde,LKH}. In addition, there has been a wealth of theoretical discoveries made relating to necessary and sufficient conditions for Hamiltonicity \cite{ore,bondy,rahman}, bounds or exact numbers of HCs in certain families of graphs \cite{schwenk,gebauer}, and the smallest non-Hamiltonian graphs satisfying certain properties \cite{holton,petersen}. HCP is also the subject of several open conjectures \cite{barnette,minreg,bridge,sheehan2,rapaport}.

However, despite this vast body of research, there has been almost no attempt made to produce a good set of HCP instances. Such sets are very important for benchmarking purposes, particularly as it is well-known that HCP can be solved in almost linear time for randomly-generated graphs \cite{frieze}. Indeed, there are good benchmark sets of instances for other important NP-complete problems such as TSP \cite{TSPLIB}, boolean satisfiability \cite{SATLIB}, and integer programming \cite{cuter}.

There are a few sets of HCP instances that exist in literature and online. The TSPLIB website \cite{TSPLIB} has a set of nine Hamiltonian instances containing between 1000 and 5000 vertices. The Extended Foster Census \cite{foster} contains a set of all 332 graphs, comprising of all known cubic symmetric graphs with less than 1000 vertices. All but five of these instances are Hamiltonian. A page hosted by the University of Kentucky \cite{kentucky} also contains eight instances of HCP with 200 vertices each, four of which are Hamiltonian. However, none of the instances in these various sets are designed to be difficult to solve, and indeed, the most effective algorithms such as Concorde \cite{concorde}, LKH \cite{LKH} and Snakes and Ladders Heuristic \cite{SLH} (SLH) are able to solve each of the instances in no more than a few seconds.

In 2015, the FHCP Challenge Set was produced to fill this void, containing 1001 instances of HCP. Since the intention was for the instances to be difficult due to their structure, rather than simply their size, no instances with more than 10000 vertices were included in the set. Every one of the 1001 instances are known to be Hamiltonian. At the 59th Annual Meeting of the Australian Mathematical Society, the {\em FHCP Challenge} was announced, in which researchers were given one year to try and solve each of the 1001 instances. A \$1001US prize was offered for the first researcher or research team to solve all of the instances, or for the reseacher who could solve the most instances in the case that nobody was able to solve all of them within a year. After twelve months, the prize was won by a two-person team consisting of David Coudert and Nathann Cohen (INRIA, France) who were successful in solving 985 instances.

Although the FHCP Challenge has now concluded, the instances contained in the FHCP Challenge Set are still of interest to researchers who want to benchmark their own HCP heuristics. In this note, a summary of the graphs contained in the FHCP Challenge Set is given, along with a short discussion about which instances proved the most difficult to solve. The FHCP Challenge Set itself, along with a solution for each instance, is available for download at http://fhcp.edu.au/fhcpcs.

\section{FHCP Challenge Set}

Graphs in the FHCP Challenge Set can be, roughly, divided into four categories, which are listed in detail below. For the most part, their difficulty comes from having relatively few Hamiltonian cycles compared to the size of the graph. However, simply having few Hamiltonian cycles is not sufficient to provide a difficult instance. For example, an arbitrarily large instance of HCP with a single Hamiltonian cycle can be easily produced by starting with a triangle, and iteratively replacing any vertex with a triangle until the desired size is reached. However, good heuristics will be able to solve instances of this type easily, because many edges can be immediately disregarded as \lq\lq bad" choices, that is, choices which cannot lead to a HC. Difficult instances, then, must be those which have structural properties that prevent the easy identification of bad edges, but such a definition is hard to define precisely. Hence, when building the FHCP Challenge Set, candidate instances were first tested on Concorde, LKH, and SLH to see if any of those algorithms struggled to solve the instance. If it was a difficult instance for at least two of the algorithms, it was included.

The four categories of instances are as follows. In each case, the vertices were randomly relabelled after construction to prevent easy identification.

\subsection{Graphs from literature}

There are a number of graph families in literature. In many cases, the exact numbers of Hamiltonian cycles in graphs from these families are known. However, practically no attention has been paid to how HCP algorithms perform on these instances. The FHCP Challenge Set included instances from the following families. First, it included generalized Petersen graphs \cite{watkins} $GP(n,2)$ for some values of $n = 3$ mod $6$. These instances are known to contain exactly three Hamiltonian cycles. It also included three families of {\em uniquely Hamiltonian graphs}, that is, graphs with only a single Hamiltonian cycle. The first was due to Sheehan \cite{sheehan}, which contains the maximal number of edges possible for uniquely Hamiltonian graphs. The second is due to Aldred and Thomassen \cite{aldred} which have average degree converging to 3 as the size of the graphs grow. The third was due to Fleischner \cite{fleischner} which is the only known family of uniquely Hamiltonian graphs with minimum degree 4 described in literature. There were some other families of uniquely Hamiltonian graphs considered, e.g. see \cite{entringer}, but these were disregarded as all three algorithms found these instances simple to solve even when the size was large.

\subsection{Modified graphs from literature}

The FHCP Challenge Set also contains two families of graphs which are derived from two families of graphs in literature. Specifically, they are derived from the generalized Petersen graphs \cite{watkins} $GP(n,2)$ for some values of $n = 5$ mod $6$, and some Flower Snarks \cite{flower}. Every graph in these two families are maximally non-Hamiltonian, that is, they are non-Hamiltonian but the addition of any edge results in a Hamiltonian graph. For the Flower Snarks, which are produced by taking $n$ copies of the star graph on 4 vertices and linking them together according to a particular scheme, a single edge is added between any two outer vertices on one of the star graphs. For the generalized Petersen graphs, if constructed as described in \cite{watkins}, a single edge is added between vertices $u_0$ and $v_{n-1}$. In both cases, it appears empirically as though the addition of this single edge results in a graph with an exponential (in $n$) number of Hamiltonian cycles. However, despite this abundance of Hamiltonian cycles, both Concorde and LKH struggled to solve these graphs, and so they were included.

\subsection{Graphs obtained by converting from other discrete problems}

HCP is an NP-complete problem, and hence any instance of any problem in NP can be transformed into an instance of HCP; however, the process of doing so often results in a large increase in the size of the instance. Conversions between NP-complete problems that result in only linear growth in the size of the instance has been a recent topic of research, e.g. see \cite{dewdney,creignou,karpred,setsplitting,hcp23hcp}. A large number of problems were converted to HCP and tested for difficulty of the resulting instance. Instances from the Dominating Set problem, the $n$-queens problem, and the generalized Instant Insanity problem were found to be sufficiently difficult to be included in the FHCP Challenge Set. Although these conversions have not been formally described in literature, source code for them is provided at \cite{FHCP}. Instances derived from a popular puzzle game, Unium \cite{unium}, were also included. In general, it seemed that SLH and LKH struggled to solve these instances, with Concorde successful for some but not others. Finally, instances of HCP which were recently derived from a famous combinatorial problem from campanology, specifically the search for peals in parts of bobs-only Stedman Triples \cite{bellringing}, were also included; none of the algorithms tested were able to easily solve these instances.

\subsection{Combinations of the above}

Since some instances are easily solved by one algorithm but not another, a number of instances were constructed by combining different instances together in various configurations. The intention was that an algorithm would need to be sufficiently clever to solve each of the component instances. Alternatively, a clever pre-processing algorithm may be successful in decomposing the instances and solving them individually.

\section{Results of the FHCP Challenge}

The {\em FHCP Challenge} ran from 30th September 2015 until 30th September 2016, after having been formally announced at the 59th Annual Meeting of the Australian Mathematical Society. The results were announced at the 60th Annual Meeting of the Australian Mathematical Society, including a summary of the top five submissions.

\begin{enumerate}\item[5.] M. Nurhafiz (Independent Researcher), 385 graphs solved.
\item[4.] M. Noisternig (TU Darmstadt, Germany), 464 graphs solved.
\item[3.] A. Gharbi and U. Syarif (King Saud University, Saudi Arabia), 488 graphs solved.
\item[2.] A. Johnson (IBM, United Kingdom), 614 graphs solved.
\item[1.] N. Cohen and D. Coudert (INRIA, France), 985 graphs solved.\end{enumerate}

As can be seen, no teams were successful in solving all instances of the FHCP Challenge Set, with all but two teams failing to solve even half of the instances. The sixteen graphs not solved by Cohen and Coudert were not solved in any of the submissions. These sixteen instance comprised three instances from \cite{bellringing}, namely {\em Sted4} (6930 vertices), {\em Sted5} (5544 vertices) and {\em Sted6} (4620 vertices), as well as thirteen instances made out of combinations of graphs which included at least one of {\em Sted4}, {\em Sted5} and {\em Sted6}. Of all the instances included in the FHCP Challenge Set, these instances were by far the most resistant to attack from heuristics. It should be noted that two other instances from \cite{bellringing} were included, namely {\em Sted10} (2772 vertices) and {\em Sted20} (1386 vertices). Of these, {\em Sted10} was only solved by Cohen and Coudert, and {\em Sted20} was only solved by Johnson, and Cohen and Coudert.

The strong result by Cohen and Coudert was obtained by spending a significant amount of time decomposing the instances and searching for similar structures between instances. Cohen and Coudert also correctly identified that many instances were made up out of smaller instances, and used solutions of the latter to obtain solutions for the former. Cohen and Coudert describe their techniques and approaches to successfully winning the FHCP Challenge in \cite{cohen}.

Some of the instances were also able to be efficiently attacked after some consideration. For instance, the uniquely Hamiltonian graphs by Sheehan are extremely dense, but a simple preprocessing algorithm can reduce the graph to a trivial instance. The other uniquely Hamiltonian graphs, as well as the generalized Petersen graphs $GP(n,2)$ with $n = 3$ mod $6$ can also be solved by recognising which family they belong to, and solving the graph isomorphism problem to obtain the original vertex labelling; this can be achieved in polynomial time since these instances all have very low maximum vertex degree \cite{luks}.

% alteratively, bibliographies prepared with BibTeX can be included by
% means of the following commands
%\bibliographystyle{srtnumbered}
%\bibliography{mybib}

\begin{thebibliography}{99}
\bibitem{concorde} D.L. Applegate, R.B. Bixby, V. Chav\'{a}tal, and W.J. Cook, {\em The Traveling Salesman Problem: A Computational Study}, Princeton University Press (2006).
\bibitem{SLH} P. Baniasadi, V. Ejov, J.A. Filar, M. Haythorpe, and S. Rossomakhine, Deterministic \lq\lq Snakes and Ladders" Heuristic for the Hamiltonian cycle problem, {\em Math. Program. Comput.}, 6(1):55-75, 2014.
\bibitem{barnette} D.W. Barnette, Conjecture 5, In: W.T. Tutte (ed), {\em Recent Progress in Combinatorics: Proceedings of the Third Waterloo Conference on Combinatorics}, New York: Academic Press, 1968.
\bibitem{bondy} J.A. Bondy and V. Chv\'{a}tal, A method in graph theory, {\em Discrete Math.}, 15(2):111-135, 1976.
\bibitem{chalaturnyk} A. Chalaturnyk, {\em A Fast Algorithm For Finding Hamiltonian Cycles.} Ph.D Thesis, University of Manitoba, 2008.
\bibitem{clark} L. Clark and R. Entringer, Smallest maximally nonhamiltonian graphs, {\em Period. Math. Hung.}, 14(1):57--68, 1983.
\bibitem{cohen} N. Cohen and D. Coudert, Le d\'{e}fi des 1001 graphes, {\em Interstices}, INRIA, 2017.
\bibitem{foster} M.D.E. Confer, P. Dobcs\'{a}nyi, B.D. McKay, and G. Royle, {\em The Extended Foster Census}, http://staffhome.ecm.uwa.edu.au/$\sim$gordon/remote/foster/, 2001, accessed April 6th 2017.
\bibitem{creignou} N. Creignou, The class of problems that are linearly equivalent to Satisfiability or a uniform method for proving NP-completeness, {\em Theor. Comput. Sci.}, 145:111-145, 1995.
\bibitem{dewdney} A.K. Dewdney, Linear transformations between combinatorial problems, {\em Internat. J. Comput. Math.}, (11):91--110, 1982.
\bibitem{gondzio} V. Ejov, J.A. Filar, and J. Gondzio, An Interior Point Heuristic for the Hamiltonian Cycle Problem via Markov Decision Processes, {\em J. Global Opt.}, 29(3):315--334, 2004.
\bibitem{branchbound} V. Ejov, J.A. Filar, M. Haythorpe and G.T. Nguyen, Refined MDP-based branch-and-bound algorithm for the Hamiltonian cycle problem, {\em Math. Oper. Res.}, 34(3):758--768, 2009.
\bibitem{hcp23hcp} V. Ejov, M. Haythorpe and S. Rossomakhine, A Linear-size Conversion of HCP to 3HCP, {\em Australasian Journal of Combinatorics}, 62(1):45--58, 2015.
\bibitem{entringer} R.C. Entringer and H. Swart, Spanning cycles of nearly cubic graphs, {\em J. Comb. Th.}, 29:303--309, 1980.
\bibitem{eppstein} D. Eppstein, The traveling salesman problem for cubic graphs, {\em J Graph Algorithm. App.}, 11(1):61--81, 2007.
\bibitem{ali} A. Eshragh, J.A. Filar, and M. Haythorpe, A hybrid simulation-optimization algorithm for the Hamiltonian cycle problem, {\em Ann. Oper. Res.}, 189(1):103--125, 2011.
\bibitem{fleischner} H. Fleischner, Uniquely Hamiltonian grpahs of minimum degree 4, {\em J. Graph Th.}, 75(2):167--177, 2014.
\bibitem{setsplitting} J.A. Filar, and M. Haythorpe, A Linearly-Growing Conversion from the Set Splitting Problem to the Directed Hamiltonian Cycle Problem, In: Honglei Xu, Xiangyu Wang (ed), {\em Optimization and Control Methods in Industrial Engineering and Construction}, Springer, Dordrech, 2014.
\bibitem{bridge} J.A. Filar, M. Haythorpe and G.T. Nguyen, A conjecture on the prevalence of cubic bridge graphs, {\em Discus. Math. Graph Th.}, 30(1):175--179, 2010.
\bibitem{karpred} J.A. Filar, M. Haythorpe, and R. Taylor, Linearly-growing reductions of Karp's 21 NP-complete problems, {\em Numer. Algebr. Control Opt.}, 8(1):1--16, 2018.
\bibitem{frieze} A. Frieze and S. Haber, An almost linear time algorithm for finding Hamilton cycles in sparse random graphs with minimum degree at least three, {\em Random Struct. Algor.}, 47(1):73--98, 2015.
\bibitem{gebauer} H. Gebauer, On the number of Hamilton cycles in bounded degree graphs, In: {\em Proceedings of the meeting on Analytic Algorithmics and Combinatorics}, Society for Industrial and Applied Mathematics, 2008.
\bibitem{cuter} N.I.M. Gould, D. Orban, and Ph. L. Toint, CUTEr (and SifDec): a Constrained and Unconstrained Testing Environment, revisited, {\em ACM Trans. Math. Soft.}, 29(4):373--394, 2003.
\bibitem{gazettepaper} M. Haythorpe, Finding Hamiltonian Cycles Using an Interior Point Method, {\em Austral. Math. Soc. Gaz.}, 37(3):170--179, 2010.
\bibitem{minreg} M. Haythorpe, On the Minimum Number of Hamiltonian Cycles in Regular Graphs, {\em Exp. Math.}, to appear, 2017. Published Online, DOI: 10.1080/10586458.2017.1306813.
\bibitem{FHCP} M. Haythorpe, {\em Flinders Hamiltonian Cycle Project}, http://fhcp.edu.au, 2017, accessed April 6th 2017.
\bibitem{bellringing} M. Haythorpe and A. Johnson, Changing Ringing and Hamiltonian Cycles : The Search for Erin and Stedman Triples, {\em Elec. J. Graph Th. App.}, submitted 2017. Available at https://arxiv.org/abs/1702.02623.
\bibitem{LKH} K. Helsgaun, An Effective Implementation of Lin-Kernighan Traveling Salesman Heuristic, {\em Eur. J. Oper. Res.} 126:106--130, 2000.
\bibitem{aldred} D.A. Holton and R.E.L. Aldred, Planar Graphs, Regular Graphs, Bipartite Graphs and Hamiltonicity, {\em Austral. J. Comb.}, 20:111-131, 1999.
\bibitem{holton} D.A. Holton and B.D. McKay, The smallest non-Hamiltonian 3-connected cubic planar graphs have 38 vertices, {\em J. Comb. Th. Ser. B.}, 45(3):305--319, 1988.
\bibitem{petersen} D.A. Holton and J. Sheehan, {\em The Petersen Graph}, Cambridge University Press, 1993.
\bibitem{SATLIB} H.H. Hoos, {\em SATLIB}, http://www.cs.ubc.ca/~hoos/SATLIB/benchm.html, 2011, accessed April 6th 2017.
\bibitem{flower} R. Isaacs, Infinite Families of Nontrivial Trivalent Graphs Which Are Not Tait Colorable, {\em Amer. Math. Monthly}, 82:221--239, 1975.
\bibitem{karp} R.M. Karp, {\em Reducibility among combinatorial problems}, Springer, New York, 1972.
\bibitem{lawler} E.L. Lawler, J.K. Lenstra, A.H.G. Rinnooy Kan and D.B. Shmoys, {\em The traveling salesman problem.} John Wiley \& Sons, Essex, England, 1985.
\bibitem{luks} E.M. Luks, Isomorphism of graphs of bounded valence can be tested in polynomial time, {\em J. Comput. Syst. Sci.}, 25:42--65, 1982.
\bibitem{ore} {\O}. Ore, Note on Hamilton circuits, {\em Am. Math. Mon.}, 67(1):55, 1960.
\bibitem{pettersson} V.H. Pettersson, Enumerating Hamiltonian Cycles, {\em Electron. J. Comb.}, 21(4):4--7, 2014.
\bibitem{rahman} M.S. Rahman and M. Kaykobad, On Hamiltonian cycles and Hamiltonian paths, {Inform. Process. Lett.}, 94:37--41, 2005.
\bibitem{rapaport} E. Rapaport-Strasser, Cayley color groups and Hamilton lines, {\em Scripta Math.}, 24:51--58, 1959.
\bibitem{TSPLIB} G. Reinelt, {\em TSPLIB}, http://www.iwr.uni-heidelberg.de/groups/comopt/software/TSPLIB95/, 2015, accessed April 6th 2017.
\bibitem{schwenk} A. Schwenk, Enumeration of Hamiltonian Cycles in Certain Generalized Petersen Graphs, {\em J. Comb. Th. Ser. B}, 47:53--59, 1989.
\bibitem{sheehan} J. Sheehan, Graphs with exactly one hamiltonian circuit. {\em J. Graph Th.}, 1:37--43, 1977.
\bibitem{sheehan2} J. Sheehan, The multiplicity of Hamiltonian circuits in a graph. {\em Recent advances in graph theory}, pp.477-480, 1975.
\bibitem{unium} J. Statz. {\em Unium}, http://www.kittehface.com/p/unium.html, 2015, accessed April 6th 2017.
\bibitem{watkins} M.E. Watkins, A Theorem on Tait Colorings with an Application to the Generalized Petersen Graphs, {\em J. Combin. Th.}, 6:152--164, 1969.
\bibitem{kentucky} {\em Hamiltonian-cycle problem}, http://www.cs.uky.edu/ai/benchmark-suite/hamiltonian-cycle.html, accessed April 6th 2017.

\end{thebibliography}

\end{document}